\documentclass[11pt, reqno]{article}
\pdfoutput=1

\usepackage{graphicx}
\usepackage{jheppub}
\usepackage{amssymb}
\usepackage{amsmath}
\usepackage[usenames,dvipsnames]{xcolor}
\usepackage{epsfig}
\usepackage{dcolumn}
\usepackage{tikz}
\usetikzlibrary{shapes.geometric, arrows}
\usepackage{upgreek}
\usepackage{setspace}
\usepackage{enumitem}
\usepackage{array,multirow,bigdelim,arydshln}
\usepackage{appendix}
\usepackage{xparse}
\usepackage[utf8]{inputenc}
\usepackage{hyperref}
\hypersetup{
	colorlinks,
	urlcolor=Maroon,
	linkcolor=Maroon,
	citecolor=Maroon
}

\usepackage{amsthm}
\newtheorem{thm}{Theorem}[section]
\newtheorem{prop}[thm]{Proposition}

\newtheorem{defn}[thm]{Definition}
\theoremstyle{definition}

\usepackage{mathtools}

\usepackage{float}
\restylefloat{table}

\NewDocumentCommand{\binomial}{omm}
 {%
  \genfrac(){0pt}{}{#2}{#3}%
  \IfValueT{#1}{_{\!#1}}%
 }
\NewDocumentCommand{\eulerian}{omm}
 {%
  \genfrac<>{0pt}{}{#2}{#3}%
  \IfValueT{#1}{_{\!#1}}%
 }

\def \s {\sigma}

\usepackage{underscore}

\usepackage{latexsym}
\usepackage{tikz}

\title{Diagonally Embedded Sets of ${\rm Trop}^+G(2,n)$'s in ${\rm Trop}\, G(2,n)$: Is There a Critical Value of $n$? }

\author{Freddy Cachazo}\emailAdd{fcachazo@pitp.ca}

\affiliation{Perimeter Institute for Theoretical Physics, Waterloo, ON N2L 2Y5, Canada}

\abstract{The tropical Grassmannian ${\rm Trop}\, G(2,n)$ is known to be the moduli space of unrooted metric trees with $n$ leaves. A positive part can be defined for each of the $(n-1)!/2$ possible planar orderings, $\alpha$, and agrees with the corresponding planar trees in the moduli space, ${\rm Trop}^{\alpha}G(2,n)$. Motivated by a physical application we study the way ${\rm Trop}^{\alpha}G(2,n)$ and ${\rm Trop}^{\beta}G(2,n)$ intersect in ${\rm Trop}\, G(2,n)$. We define their intersection number as the number of unrooted binary trees that belong to both and construct a $(n-1)!/2\times (n-1)!/2$ intersection matrix. We are interested in finding the diagonal (up to permutations of rows and columns) submatrices of maximum possible rank for a given $n$. We prove that such diagonal matrices cannot have rank larger than $(n-3)!$ using the CHY formalism. We also prove that the bound is saturated for $n=5$ (the condition is trivial for $n=4$), that for $n=6$ the maximum rank is $4$, and that for $n=7$ the maximum rank is $\geq 14$. We also ask the following question: Is there a value $n_{\rm c}$ so that for any $n>n_{\rm c}$ the bound $(n-3)!$ is always saturated? We review and extend two relevant results in the literature. The first is the Kawai-Lewellen-Tye (KLT) choice of sets which leads to a $(n-3)!\times (n-3)!$ block diagonal submatrix with blocks of size $d\times d$ with $d = \lceil (n-3)/2\rceil !\lfloor (n-3)/2\rfloor !$. The second result is that the number of ${\rm Trop}^\alpha G(2,n)$'s that intersect a given one grows as $\exp(n\log (3+\sqrt{8}))$ for large $n$ which implies that the density of the intersection matrix goes as $\exp(-n(\log(n)-2.76))$. We interpret this as an indication that the generic behavior is not seen until $n \approx \exp (2.76)$, i.e. $n = 16$. We also find an exact formula for the number of zeros in a KLT block.}

\begin{document}
\maketitle
\addtocontents{toc}{\protect\setcounter{tocdepth}{1}}
\def \tr {\nonumber\\}
\def \nn {\nonumber}
\def \la {|}
\def \ra {|}
\def \dd {\Theta}
\def\hset{\texttt{h}}
\def\gset{\texttt{g}}
\def\sset{\texttt{s}}
\def\A {\textsf{A}}
\def\B {\textsf{B}}
\def\C {\textsf{C}}
\def\D {\textsf{D}}
\def\E {\textsf{E}}
\def\F {\textsf{F}}
\def\G {\textsf{G}}
\def\I {\textsf{I}}
\def\J {\textsf{J}}
\def\H {\textsf{H}}
\def \be {\begin{equation}}
\def \ee {\end{equation}}
\def \ba {\begin{eqnarray}}
\def \ea {\end{eqnarray}}
\def \k {\kappa}
\def \h {\hbar}
\def \r {\rho}
\def \l {\lambda}
\def \be {\begin{equation}}
\def \en {\end{equation}}
\def \bes {\begin{eqnarray}}
\def \ens {\end{eqnarray}}
\def \red {\color{Maroon}}
\def \pt {{\rm PT}}
\def \s {\textsf{s}}
\def \t {\textsf{t}}
\def \C {\textsf{C}}
\def \tp {||}
\def \p {x}
\def \x {z}
\def \V {\textsf{V}}
\def \ls {{\rm LS}}
\def \ma {\Upsilon}
\def \SL {{\rm SL}}
\def \GL {{\rm GL}}
\def \w {\omega}
\def \e {\epsilon}
\def \a {\alpha}
\def \b {\beta}
\def \ort {\textsf{O}}

\numberwithin{equation}{section}

\section{Introduction: The Question}

Tropical Grassmannians, ${\rm Trop}\, G(k,n)$, were introduced by Speyer and Sturmfelds in \cite{SSTrop} who also showed that ${\rm Trop}\, G(2,n)$ is the moduli space of phylogenetic trees studied by Billera, Holmes and Vogtmann (BHV) \cite{BHV}. Motivated by the work of Postnikov \cite{Alex} on totally positive Grassmannians, $G^+(k,n)$, Speyer and Williams introduced positive tropical Grassmannians  ${\rm Trop}^+\, G(k,n)$ \cite{SWTrop}. Moreover, they also showed that ${\rm Trop}^+\, G(2,n)$ is closely related to the associahedron and it is the moduli space of planar trees.

In this work we are only interested in ${\rm Trop}\, G(2,n)$. It is well-known that in order to construct a positive part of ${\rm Trop}\, G(2,n)$ one has to select an ordering $\alpha$ of the set of $n$ elements $\{1,2,\ldots n\}$. We introduce the notation ${\rm Trop}^{\alpha}\, G(2,n)$ to refer to the corresponding space of planar trees. Any two such orderings related by a cyclic transformation give rise to the same object in ${\rm Trop}\, G(2,n)$, so we take $\a := (\a_1,\a_2,\ldots ,\a_n )$ to be an element of $S_n/\mathbb{Z}_n$ which defines an ordering.  Moreover we will often abuse notation and refer to ${\rm Trop}^{\alpha}\, G(2,n)$ by its ordering $\alpha$.

The connection between ${\rm Trop}\, G(2,n)$ and the BHV space of phylogenetic trees is very useful. In a nutshell, the BHV construction assigns an $(n-3)$-dimensional orthant to each unrooted binary tree with $n$ leaves. Each direction corresponds to the length of an edge in the tree. Codimension one boundaries correspond to trees with one zero-length edge, i.e. with a vertex of degree four. There are three ways of opening up (or smoothing) a degree four vertex into two degree three vertices and an edge. One of them leads back the original orthant and the other two connect to two other orthants. The total number of orthants is the number of binary trees $(2n-5)!!$. For each orthant we can assign coordinates $x_i$ and consider the hypersurface $\sum_i x_i =1$. The union of all such hypersurfaces in the BHV space is known as the {\it link of the origin} \cite{BHV}.

Choosing an ordering $\a = (\a_1,\a_2,\ldots ,\a_n )$ means restricting to planar trees with respect to the ordering. There are  $C_{n-2}$ planar binary trees, with $C_m$ the $m^{\rm th}$ Catalan number. The link of the origin of the corresponding orthants is the boundary of the dual associahedron with $n$ letters (See Prop. 3.1 in \cite{BHV}). Using the identification with ${\rm Trop}^{\alpha}\, G(2,n)$, this construction gives an embedding of it into the full ${\rm Trop}\, G(2,n)$.   

There are clearly $(n-1)!/2$ distinct ${\rm Trop}^\alpha\, G(2,n)$ in ${\rm Trop}\, G(2,n)$. Note that an ordering and its reflection are equivalent and hence it is necessary to divide by two. The pattern of how these spaces of planar trees intersect in ${\rm Trop}\, G(2,n)$ is the main topic of this work and can be formalized by defining the intersection matrix.

\begin{defn}\label{inter}
    Given ${\rm Trop}^\alpha\, G(2,n)$ and ${\rm Trop}^\beta\, G(2,n)$ defined by their orderings $\a$ and $\b$, their intersection number,  $I(\a,\b)$, is the number of unrooted binary trees that belong to both. We call the $(n-1)!/2\times (n-1)!/2$ matrix with entries $I(\a,\b)$ the intersection matrix. 
\end{defn}

For example, if there is not any unrooted binary tree which is planar with respect to both orderings then $I(\a,\b)=0$. The self-intersection is $I(\a,\a)=C_{n-2}$ and it is clearly the maximum possible value of $I(\a,\b)$. 

In this work we are mainly interested in certain submatrices of the intersection matrix defined as follows: 

\begin{defn}\label{diagonal}
    A diagonal submatrix of the intersection matrix is given by two sets of orderings ${\cal A}$ and ${\cal B}$, not necessarily the same, with $|{\cal A}|=|{\cal B}|$ and which satisfy that for $\a_\I \in {\cal A}$ and $\b_\J \in {\cal B}$ with $\I,\J\in \{1,2,\ldots ,|{\cal A}|\}$, the matrix $I(\a_\I,\b_\J)$ is diagonal, up to permutations of rows and columns, and has rank $|{\cal A}|$.  
\end{defn}

In section 2 we prove the following proposition.

\begin{prop}\label{rank}
  The rank of a diagonal submatrix of the intersection matrix in ${\rm Trop}\, G(2,n)$ cannot be larger than $(n-3)!$.  
\end{prop}

The proof of Proposition \ref{rank} uses what is known as the Cachazo-He-Yuan (CHY) construction \cite{Cachazo:2013gna,Cachazo:2013hca,Cachazo:2013iea} which connects the moduli space of punctured Riemann spheres with unrooted binary trees. Moreover, the CHY construction also gives a very clear interpretation of the Kawai-Lewellen-Tye (KLT) relation between quantum field theories such as Einstein gravity and Yang-Mills theory \cite{KLT}. One of the motivations for this work is the fact that the KLT relations would dramatically simplify if a diagonal submatrix of rank $(n-3)!$ of the intersection matrix existed. 

This led us to start the study of diagonal submatrices. Let us emphasize again that the rows and columns of the submatrix do not have to correspond to the same sets of orderings. Knowing that the rank is bounded by $(n-3)!$, it is natural to ask if the bound is always saturated. We show in section 2 that it is saturated for $n=5$ (the case $n=4$ is trivial) but it is not for $n=6$. In fact, we prove that the maximum rank for $n=6$ is $4$. This motivates the following definition.   

\begin{defn}\label{troprank}
    The diagonal degree of ${\rm Trop}\, G(2,n)$ is the maximum rank of all possible diagonal submatrices of the intersection matrix $I(\a,\b)$. 
\end{defn}

For example, using this terminology the diagonal degree of ${\rm Trop}\, G(2,6)$ is $4$. In section 5 we show that $14$ is a lower bound for the diagonal degree of ${\rm Trop}\, G(2,7)$.

The main question we pose in this note is the following.

\vspace{0.1in}

\noindent {\bf Question 1.4.} {\it Is there a value $n_{\rm c}$ such that for any $n>n_{\rm c}$ the diagonal degree of ${\rm Trop}\, G(2,n)$ is $(n-3)!$, i.e. such that the bound is saturated?}

\vspace{0.1in}

If such a value $n_c$ exists we call it the {\it critical value}. At this point we do not have much evidence to support the existence of $n_c$. Instead, in sections 3 and 4 we review and extend some facts known in the literature which indirectly hint to a positive answer to Question 1.4. 

Let us illustrate the definitions with an example.

\vspace{0.1in}

\noindent {\bf Example 1.5.} The case $n=5$ is the only one where the bound is known to be saturated ($n=4$ is trivial). Let us prove this by explicitly constructing a $2\times 2$ diagonal submatrix of the $12\times 12$ intersection matrix. The simplest way to select the two sets ${\cal A}$ and ${\cal B}$, each with $(5-3)!=2$ orderings is 
\be
{\cal A} = {\cal B} = \{ (1,2,3,4,5), (1,3,5,2,4) \}.
\ee
The corresponding intersection submatrix is easily computed and shown to be a $2\times 2$ diagonal matrix by noting that 
\be\label{ortho}
I((1,2,3,4,5), (1,3,5,2,4)) = 0
\ee
and
\be
I((1,2,3,4,5),(1,2,3,4,5))=I((1,3,5,2,4),(1,3,5,2,4))=5.
\ee
The explicit form of the diagonal submatrix is 
\be
\left( \begin{array}{cc}
    5 & 0  \\
    0 & 5  \\
  \end{array} \right).
\ee

One way to see why \eqref{ortho} is true is by explicilty showing both ${\rm Trop}^{(1,2,3,4,5)}G(2,5)$ and ${\rm Trop}^{(1,3,5,2,4)}G(2,5)$ in ${\rm Trop}\, G(2,5)$. In figure 1 we show the link of the origin in the BHV space of trees which happens to be a Petersen graph in this case. Positive tropical Grassmannians are cycles of length five embedded in the Petersen graph.

    \begin{figure}[h!]
	\hspace{0.1in}
    \includegraphics[width=0.9\linewidth]{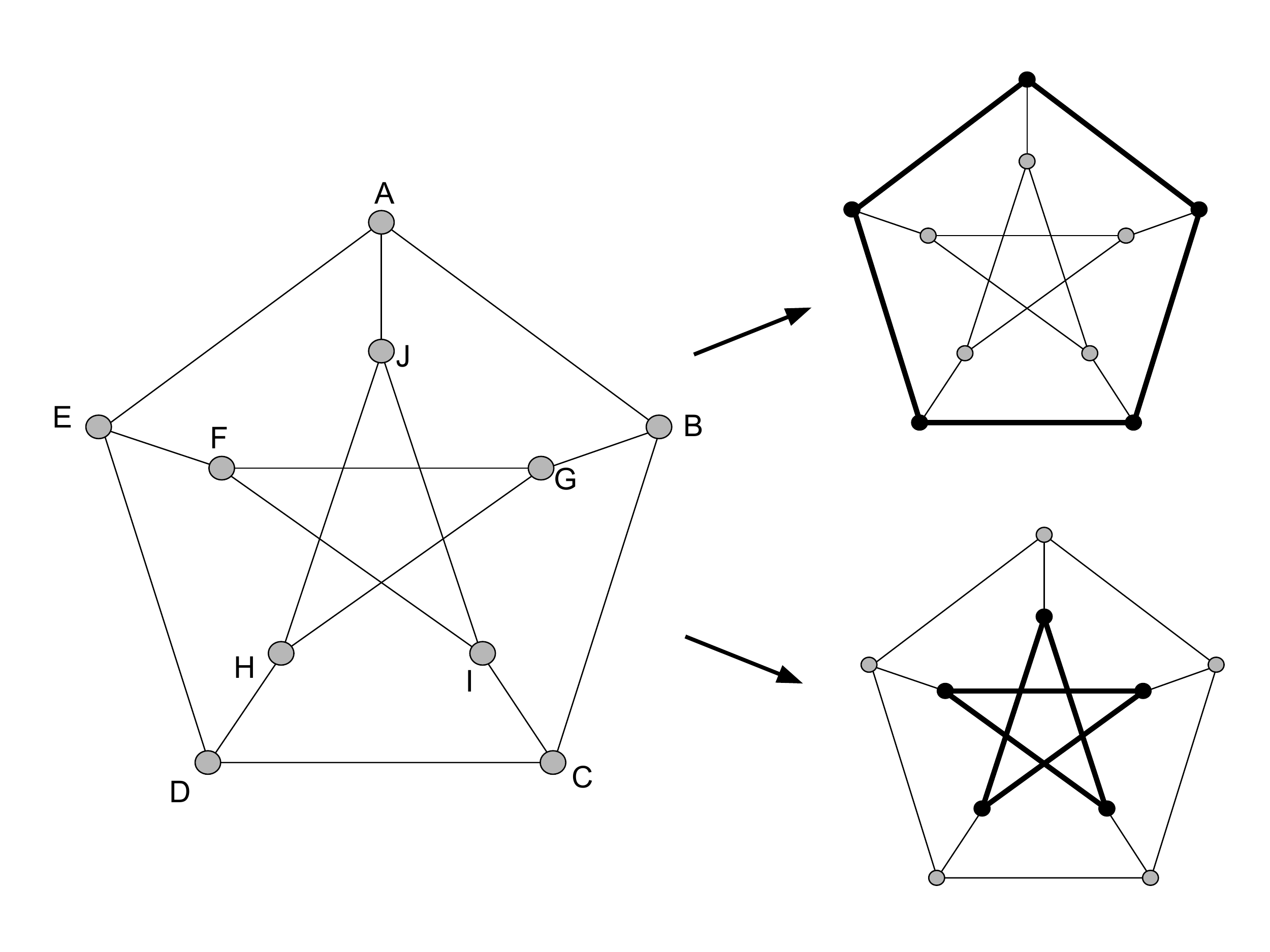}
	\caption{The intersection of ${\rm Trop}\, G(2,5)$ with a unit sphere is isomorphic to the Petersen graph. This is also the link of the origin in the BHV space of trees. On the left, vertices have been labeled according to the tree they represent: $\A\rightarrow (34)(512)$, $\B\rightarrow (12)(345)$, $\C\rightarrow (45)(123)$, $\D\rightarrow (23)(451)$, $\E\rightarrow (51)(234)$, $\F\rightarrow (24)(135)$, $\G\rightarrow (35)(241)$, $\H\rightarrow (41)(352)$, $\I\rightarrow (13)(524)$, and $\J\rightarrow (52)(413)$. The notation is explained in the text. On the right, two different ${\rm Trop}^+G(2,5)$ are shown, corresponding to orderings $(1,2,3,4,5)$ and $(1,3,5,2,4)$. Note that they do not intersect. }
	\label{fiveNaive}
    \end{figure}

Each edge in figure 1 represents one unrooted binary tree with five leaves. Each vertex in the graph is a tree with one degree three and one degree four vertex. We use the notation $(ab)(cde)$ for such trees, where $(ab)$ are the leaves of the degree three vertex and $(cde)$ those of the degree four vertex. While the figure only labels the trees on the vertices, this is enough to uniquely fix the binary trees on the edges. For example, consider the edge $\overline{\A\B}$, the vertices are $(34)(512)$ and $(12)(345)$. The only tree that can degenerate to both is $(12)(34)5$, where the grouping indicates the leaves attached to the same trivalent vertex. The embedding of the ${\rm Trop}^+G(2,5)$'s given by $(1,2,3,4,5)$ and $(1,3,5,2,4)$ is easily found by cycling their labels to make five combinations of the form $(ab)(cde)$ and finding the corresponding vertices.   

While this is the simplest way to get a maximal diagonal submatrix for $n=5$, there is an alternative one known as the KLT choice in which the two sets, ${\cal A}$ and ${\cal B}$, do not share any element. In fact, in physical applications the KLT choice leads to important simplifications in calculations. The reason is that in the KLT choice whenever a ${\rm Trop}^+G(2,5)$'s in ${\cal A}$ has a non-zero intersection with one in ${\cal B}$, the intersection number is exactly $1$, i.e. they share a single unrooted binary tree. We discuss the KLT choice in detail in section 3.

\vspace{0.1in}

This note is organized as follows: In section 2 we prove that the rank of a diagonal submatrix of the intersection matrix cannot be larger than $(n-3)!$ and study the $n=6$ case where the bound is not saturated and the maximum rank is $4$. In sections 3 and 4 we review some relevant results already known in the literature which show that the existence of a critical value is plausible. More explicitly, in section 3 we review the KLT construction which leads to block diagonal matrices of size $(n-3)!\times (n-3)!$ with blocks of size $d\times d$, where $d = \lceil (n-3)/2\rceil !\lfloor (n-3)/2\rfloor ! $. In section 4 we show that the density of the intersection matrix goes as ${\rm exp}(-n(\log n -2.76))$ for large $n$ and also provide an explicit formula for the density of KLT blocks for all $n$. In section 5 we discuss some future directions and provide tools such as the CHY graphical technique for counting intersections and an algorithm for finding permutation submatrices of a binary matrix. In appendix A we review the elements of the CHY formalism needed for the proof of Proposition \ref{rank} and in appendix B we provide a derivation of the asymptotic behavior of super Catalan numbers.    

\section{Upper Bound on the Rank of Diagonal Submatrices}

In this section we prove Proposition \ref{rank} which states that the intersection matrix, introduced in Definition \ref{inter}, cannot have a diagonal submatrix of rank larger than $(n-3)!$. 

The proof uses the CHY formulation of scattering amplitudes. A short review with the relevant constructions is given in appendix A. Here all we need is that the CHY formulation associates a real $(n-3)!$-dimensional vector to each ordering $\alpha\in S_n/\mathbb{Z}_n$. These $(n-1)!/2$ vectors have the property that their inner product computes a rational function, known in physics as a scattering amplitude, which vanishes if and only if there are no unrooted binary trees with $n$ leaves which are planar with respect to both orderings. It then follows that letting $\phi(\alpha)$ and $\phi(\beta)$ be the vectors associated with $\alpha,\beta \in S_n/\mathbb{Z}_n$ then $\phi(\a)\cdot \phi(\b) = 0$ if and only if $I(\a,\b)=0$.

Now we are ready to prove Proposition \ref{rank}.

\begin{proof}
  Let us proceed by contradiction. Assume that there is a diagonal submatrix of the intersection matrix $I(\a,\b)$ with rank $r>(n-3)!$. This implies that the $(n-1)!/2\times (n-1)!/2$ matrix $m(\a,\b):=\phi(\a)\cdot \phi(\b)$ also possesses a diagonal submatrix of rank $r>(n-3)!$. However, this is impossible since the rank of $m(\a,\b)$ is at most $(n-3)!$. This is because $m(\a,\b)$ is the Gram matrix of vectors in $\mathbb{R}^{(n-3)!}$.  
\end{proof} 

Before ending this section let us show that the bound $(n-3)!$ is not saturated for $n=6$. 

\vspace{0.1in}

\noindent {\bf Example 2.1.} We studied ${\rm Trop}\, G(2,6)$ and constructed the corresponding $60\times 60$ intersection matrix. Performing an exhaustive computer-assisted search for diagonal submatrices is not difficult\footnote{In section 5 we explain the algorithm used to carry out this search.}. We found that there are not any diagonal submatrices of rank $5$ and found an example of a diagonal submatrix of rank $4$. Let the two sets of ${\rm Trop}^+ G(2,6)$'s be 
\begin{align}
    {\cal A} & =  \{ (1,2,3,4,5,6),(1,2,6,3,4,5),(1,5,3,4,2,6),(1,5,4,2,3,6) \},  \\
    {\cal B} & =  \{ (1,3,2,5,6,4),(1,2,5,3,6,4),(1,3,5,2,6,4),(1,3,6,5,2,4) \}.
\end{align}
The corresponding intersection matrix is 
\be
\left( \begin{array}{cccc}
    1 & 0 & 0 & 0  \\
    0 & 1 & 0 & 0  \\
    0 & 0 & 1 & 0  \\
    0 & 0 & 0 & 1
  \end{array} \right).
\ee

\vspace{0.1in}

Having shown why diagonal submatrices of rank $(n-3)!$ of the intersection matrix are special if they exist, we now provide some facts that might be useful in finding them.

\section{KLT Block Diagonal Sets}

In this section we present a choice of sets that gives rise to a block diagonal submatrix of the intersection matrix of size $(n-3)!\times (n-3)!$. The blocks are of size $d\times d$ with $d = \lceil (n-3)/2\rceil !\lfloor (n-3)/2\rfloor !$. Note that for $n=5$ this gives rise to a diagonal matrix as we show in Example 3.2. 

The choice of sets comes from the Kawai-Lewellen-Tye (KLT) construction that has its origins in string theory and was developed in the 80's as a way of relating physical quantities of closed strings to those of open strings \cite{KLT}. The explicit form of the sets reviewed here was first presented for all values of $n$ in \cite{BernKLT}.

Consider the first set of $(n-3)!$ orderings to be
\be
{\cal A} := \{ (1,\omega(2),\omega(3),\ldots ,\omega(n-2) ,n-1,n) : \, \omega \in S_{n-3} \}
\ee
where the permutations $\omega \in S_{n-3}$ act on the set $\{ 2,3,\ldots ,n-2\}$. The second set is  
\be
{\cal B} := \{ (1,\gamma(2),\ldots \gamma(m), n,\gamma(m+1),\ldots \gamma(n-2),n-1) : \, \gamma \in S_{n-3} \}
\ee
with $m = \lceil (n-3)/2\rceil +1$ and $\gamma \in S_{n-3}$ acting on the set $\{ 2,3,\ldots ,n-2\}$. 

\begin{prop}\label{bloki}
  ${\rm Trop}^{\a}G(2,n)$ with $\a =(1,\omega(2),\omega(3),\ldots ,\omega(n-2) ,n-1,n) \in {\cal A}$ has vanishing intersection with ${\rm Trop}^{\b}G(2,n)$ with $\b =(1,\gamma(1),\ldots \gamma(m), n,\gamma(m+1),\ldots \gamma(n-2),n-1)\in {\cal B}$ if 
  \be
  \{ \omega(2),\omega(3),\ldots ,\omega(m)\} \cap \{ \gamma(m+1),\gamma(m+2),\ldots \gamma(n-2)\} \neq \emptyset .
  \ee
\end{prop}

This fact is an implicit consequence of the string theory construction using vertex operators. However, there is a simple combinatorial proof which we provide for the reader's convenience.

\begin{proof}
  Let us consider first the cases with $n>5$ and take any orderings where there is a non-empty intersection, i.e. there exists 
  \be
  \rho\in \{ \omega(2),\omega(3),\ldots ,\omega(m)\} \cap \{ \gamma(m+1),\gamma(m+2),\ldots \gamma(n-2)\}.
  \ee
  Now proceed by contradiction assuming that there exists at least one binary tree, ${\cal T}$, which is planar with respect to both orderings. Note that in ${\cal T}$ no pair of labels in $\{n-1,n,1\}$ are adjacent in both orderings so they cannot be neighbors, i.e. no pair can belong to the same cherry (a trivalent vertex with two external edges). It is also easy to see that in fact $n$ cannot be in any cherry. This means that ${\cal T}$ must contain at least one cherry with either no elements in $\{n-1,1\}$ or at most one. Let use denote the leaves in such a cherry $\{ e_1,e_2\}$. Therefore these labels are adjacent (either as $(\ldots,e_1,e_2,\ldots )$ or $(\ldots,e_2,e_1,\ldots )$ ) in both orderings. 
  
  Removing (or pruning) the cherry with leaves $\{ e_1,e_2\}$ in ${\cal T}$ gives rise to a tree ${\cal T}'$ with $n-1$ leaves. The new tree has $n-2$ leaves inherited from ${\cal T}$ and a new leaf we call $e_{\rm new}$. If either $1$ or $n-1$ were part of the cherry then we set $e_{\rm new}$ to that value. The new tree is planar with respect to the two orderings obtained by replacing $(\ldots,e_1,e_2,\ldots )$ and $(\ldots,e_2,e_1,\ldots )$  by $(\ldots ,e_{\rm new},\ldots )$. The new pair of orderings again share the same property as the original one, i.e. there exists a
  \be
  \rho'\in \{ \omega'(2),\omega'(3),\ldots ,\omega'(m')\} \cap \{ \gamma'(m'+1),\gamma'(m'+2),\ldots \gamma'(n-3)\}.
  \ee  
  Repeating this procedure to produce trees with fewer and fewer leaves one arrives at a tree with five leaves. This is the point where the condition of having a non-empty intersection is crucial as we arrive at two orderings $(1,a,b,n-1,n)$ and $(1,c,n,d,n-1)$ satisfying $\{a\} \cap \{d\} \neq \emptyset$, i.e. $a=d$. This implies that $c=b$. Rewriting the two orderings gives $(1,a,b,n-1,n)$ and $(1,b,n,a,n-1)$, but we have seen in figure 1 that pairs like this have zero intersection in the BHV space of trees with five leaves. This means that there is no tree which is planar with respect to both ordering and hence we have reached a contraction implying that the original ${\cal T}$ did not exist. 
\end{proof}

Proposition \ref{bloki} motivates the separation of the $(n-3)!$ orderings in ${\cal A}$ and ${\cal B}$ each into $(n-3)!/d$ sets with $d = \lceil (n-3)/2\rceil !\lfloor (n-3)/2\rfloor !$. The construction starts by partitioning $\{2,3,\ldots ,n-2\}$ into two sets of sizes $\lceil (n-3)/2\rceil !$ and  $\lfloor (n-3)/2\rfloor !$ in all possible ways. Let us denote the different sets as $H_\texttt{I}$ and $G_\texttt{I}$, i.e. $H_\texttt{I}\,\cup\, G_\texttt{I} = \{2,3,\ldots ,n-2\}$ and $\texttt{I}\in \{ 1, \ldots ,(n-3)!/d\}$. Now, let ${\cal A}_\texttt{I}\subset {\cal A}$ be given by orderings of the form $(1,H_\texttt{I},G_\texttt{I},n-1,n)$ and ${\cal B}_\texttt{I}\subset {\cal B}$ be given by orderings of the form $(1,H_\texttt{I},n,G_\texttt{I},n-1)$. Using Proposition \ref{bloki}, it is clear that the intersection matrix of ${\cal A}$ and ${\cal B}$ is block diagonal. In other words, we have 
\be
I(\a\in {\cal A}_\texttt{I}  ,\b \in {\cal B}_\texttt{J} ) = 0 \qquad {\rm if} \quad \texttt{I}\neq \texttt{J}.
\ee
Of course, inside each of the blocks, i.e. when $\texttt{I}=\texttt{J}$, there can be zero entries but only if $n$ is large enough as we discuss in the next section.

\vskip0.1in

    \begin{figure}[h!]
    \hspace{-0.35in}
    \includegraphics[width=1.05\linewidth]{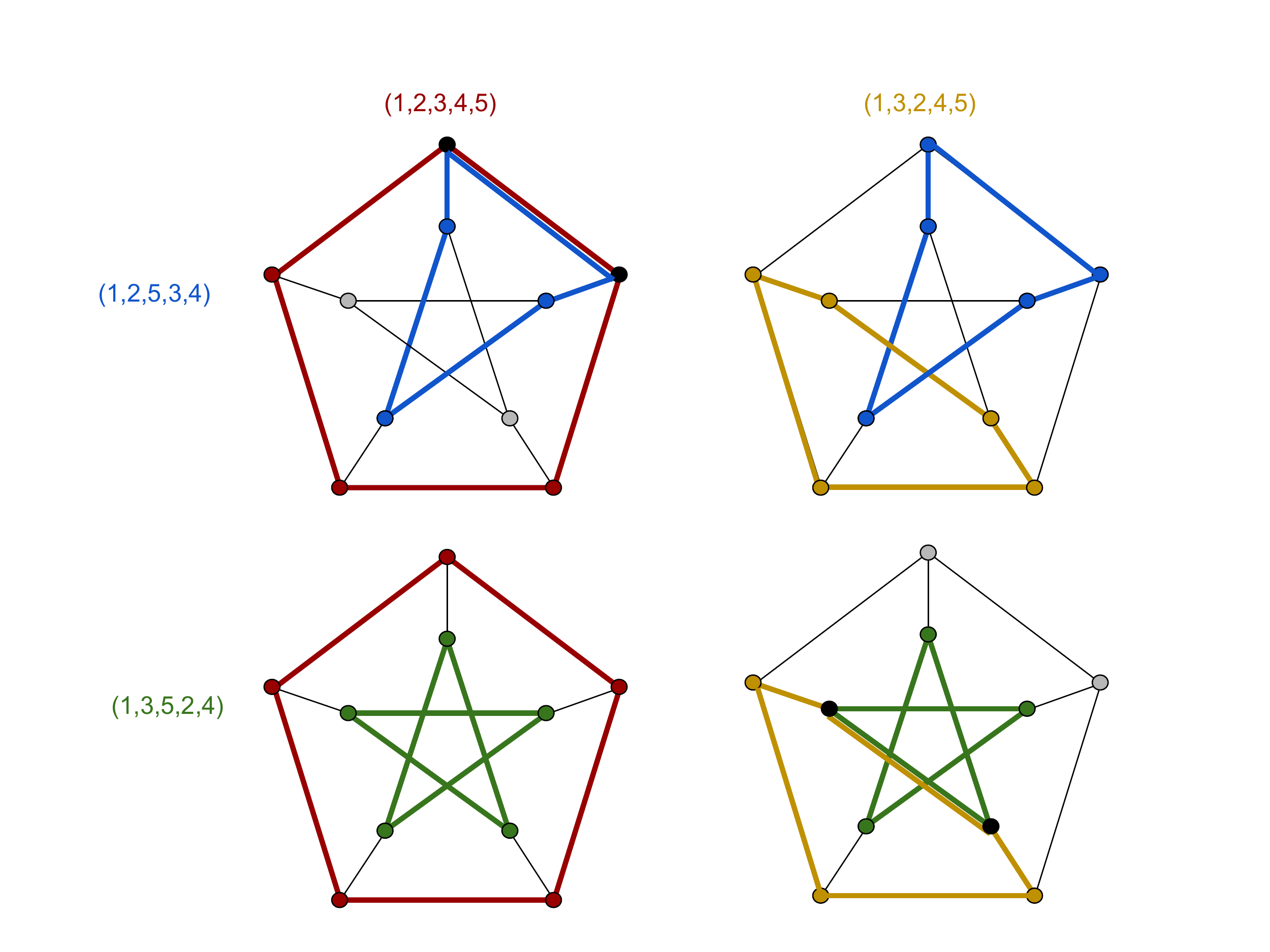}
	\caption{Each entry in this $2\times 2$ matrix shows a ${\rm Trop}\, G(2,5)$ with two ${\rm Trop}^+ G(2,5)$ embedded. The embedding has been done using the labeling from figure 1. The simplest approach is to find the five vertices that belong to a given ${\rm Trop}^+ G(2,5)$. For example, $(1,2,5,3,4)$ contains vertices $(125)(34)$, $(253)(41)$, $(534)(12)$, $(341)(25)$, and $(412)(53)$.}
	\label{fiveKLT}
    \end{figure}

\noindent {\bf Example 3.2.} The $n=5$ sets ${\cal A}$ and ${\cal B}$ given in Example 1.5 do not coincide with the KLT ones. The KLT sets are given by 
\be
{\cal A} = \{ (1,2,3,4,5), (1,3,2,4,5) \}\quad {\rm and} \quad  {\cal B} = \{ (1,2,5,3,4), (1,3,5,2,4) \}.
\ee

In figure \ref{fiveKLT} we arrange the embedding of one element in ${\cal A}$ and one from ${\cal B}$ in ${\rm Trop}\, G(2,5)$ into a $2\times 2$ matrix. The intersection matrix is clearly the identity matrix since the off-diagonal terms do not intersect while the elements in the diagonal intersect in exactly one binary tree.

\section{Sparsity of the Intersection Matrix at Large n}

In this section we discuss the density of the intersection matrix. By this we mean the ratio of non-zero entries to the total number. The first observation is that the rows and columns of the intersection matrix are permutations of each other. In particular, they all contain the same number of zeroes. This is why it is interesting to find out how many entries in a row are zeroes or equivalently how many entries are non-zero.    

This problem can be addressed by mapping it to a dual version involving associahedra as done by Mizera in \cite{Sebastian1} building on \cite{Deva1,Deva2}. Mizera notices that the real part of the moduli space of punctured Riemann spheres, $\overline{{\cal M}}_{0,n}(\mathbb{R})$ can be tiled by $(n-1)!/2$ associahedra and finds that the number of associahedra that intersects a single associahedron in codimension $k$ facets is given by $T(n-1,k+1)$, where $T(m,r)$ is the number of diagonal dissections of a convex $m$-gon into $r+1$ regions (see e.g. the sequence A033232 in the OEIS \cite{OEIS}). This means that in order to find the total number of associahedra intersecting a given one it is enough to sum $T(n-1,k+1)$ over all values of $k$. The number of all such subdivisions is well-known and it is given by the super Catalan or Schröder–Hipparchus numbers $\texttt{S}_{n-1}$. The first few corresponding to $n=4,5,6,7,8$ are $3, 11, 45, 197, 903, 4279$ respectively (see e.g. the sequence A001003 in the OEIS \cite{OEIS}). 

In \cite{CYY}, Yeats, Yusim, and the author used the fact that the number of non-zero entries is given by  $\texttt{S}_{n-1}$ in order to study the density of a matrix intimately related to the intersection matrix. In fact, their result can be directly used in our context. 

Let us see how $\texttt{S}_{n-1}$ compares to $(n-1)!/2$, the total number of ${\rm Trop}^+ G(2,n)$ inside ${\rm Trop}\, G(2,n)$. This is relevant to Question 1.4 in the asymptotic regime when $n$ is large. The asymptotic behavior of the super Catalan numbers is known to be \cite{OEIS}
\be
\log(\texttt{S}_{n-1})\sim n\, \log\left(3+\sqrt{8}\right)-\frac{3}{2}\log(n) + {\cal O}(n^0). 
\ee
This results is also derived in Appendix B for the reader's convenience. 

This number is very small compared to the total number of orderings
\be
\log((n-1)!/2)\sim n\,\log (n)-n-(1/2)\log(n)+{\cal O}(n^0). 
\ee
In other words, this shows that asymptotically the intersection matrix's density is given by 
\be\label{asym}
\frac{\texttt{S}_{n-1}}{(n-1)!/2} \sim e^{-n(\log (n) - 2.76)}.
\ee

The sparsity of a matrix is usually defined to be the ratio of the number of zero entries to the total number of entries or in this case $1$ minus the density. This means that
\be
{\rm sparsity}\, \sim 1- e^{-n(\log (n) - 2.76)}
\ee
which goes to $1$ for large $n$ and therefore we can say that the intersection matrix is asymptotically sparse.

While the asymptotic behavior is clear, it is interesting to see how it is reached. This is possible since we have exact formulas for \cite{wolfram}
$$ \texttt{S}_{n-1}  = \frac{3P_{n-2}(3)-P_{n-3}(3)}{4(n-1)}, $$
where $P_m(x)$ are Legendre polynomials.

\begin{figure}[h!]
	\hspace{0.35in}
    \includegraphics[width=0.8\linewidth]{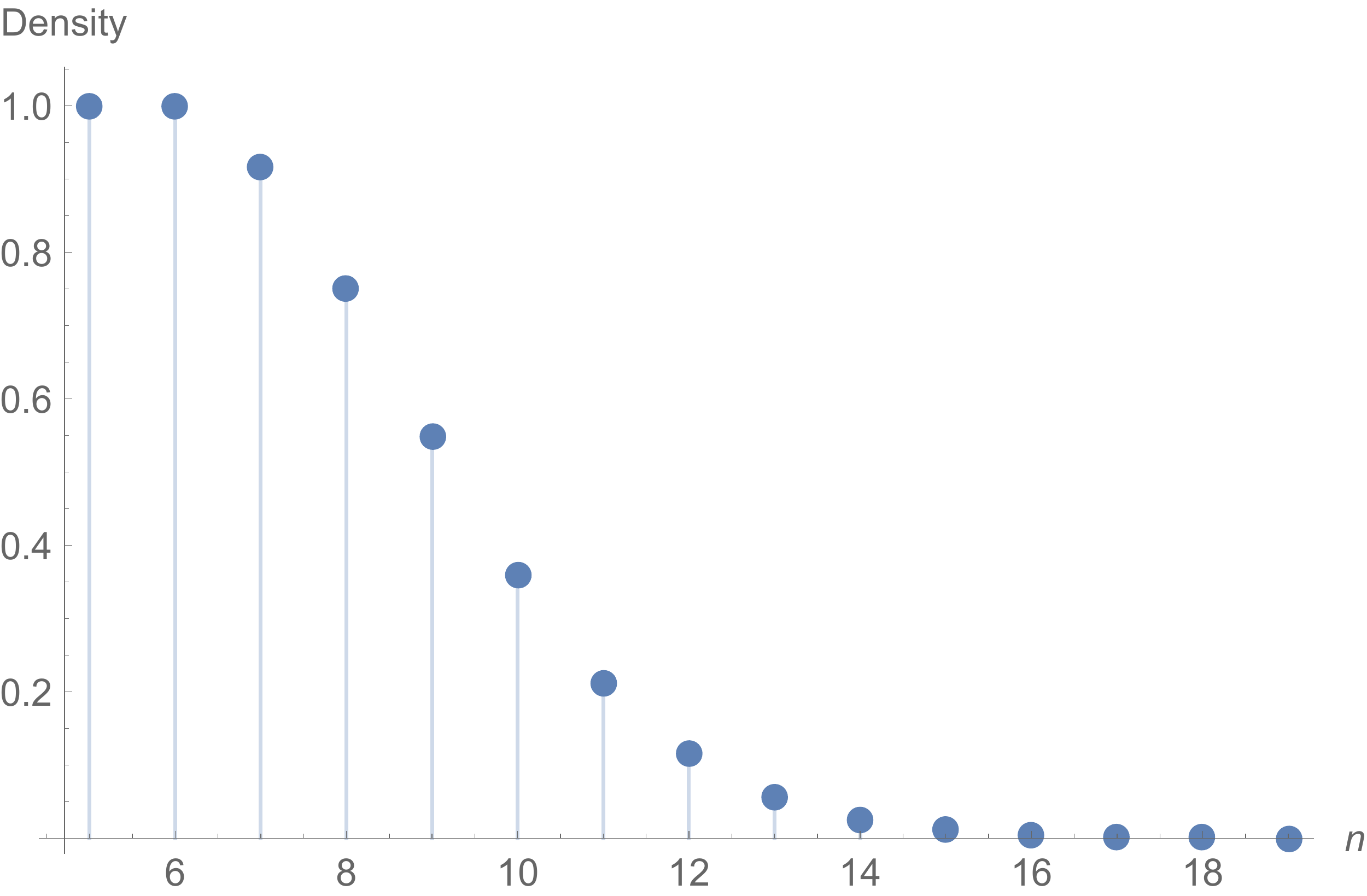}
	\caption{This plot shows the ratio of the number of non-zero entries in the intersection matrix to the total number of entries. This shows that the generic behaviour is only reached for values of $n$ larger than $16$.}
	\label{density}
    \end{figure}

The plot of the ratio $\texttt{S}_{n-1}/((n-1)!/2)$ for $n=5,6,\ldots ,19$ is shown in figure \ref{density} and seems to indicate that the generic behaviour is reached for $n>16$. 

\subsection{Behaviour of KLT Blocks at Large $n$}

Let us end this section by revisiting the KLT sets from section 3 and performing an asymptotic analysis. Recall that the KLT sets give rise to a $(n-3)!\times (n-3)!$ submatrix of the intersection matrix which is block diagonal. Here we are interested in the structure of a single block. Each block is a $d\times d$ matrix with $d = \lceil (n-3)/2\rceil !\lfloor (n-3)/2\rfloor !$.

All blocks are permutations of each other so it is sufficient to consider the block containing the canonical ordering. Once again we are interested in the density of the block, i.e., the ratio of non-zero entries to the total number of entries. For the reader familiar with the KLT literature it might be surprising that a KLT block can have zeros in it. Recall that in physical applications the intersection matrix has to be replaced by the matrix of biadjoint scalar amplitudes $m(\alpha,\beta)$. While the KLT blocks in this matrix of amplitudes can have zeroes, the object that enters in most computations is the inverse matrix of the blocks and those are known to have no zeroes in them for any values of $n$. Of course, this is not a contradiction since it is well-known that sparsity of a matrix is a property that is not generically preserved by the inverse matrix operation.

Within a KLT block, each row is a permutation of the first one so we can once again count the number of intersections between $(1,2,\ldots ,m,m+1,\ldots , ,n)$ and orderings of the form
\be
(1,\gamma(2),\ldots \gamma(m), n,\rho(m+1),\ldots \rho(n-2),n)
\ee
with $\gamma$ permutations of the set $\{ 2,3,\ldots ,m\}$, $\rho$ permutations of the set $\{ m+1,m+2,\ldots ,n-2\}$ and $m = \lceil (n-3)/2\rceil +1$.

    \begin{figure}[h!]
	\hspace{0.35in}
    \includegraphics[width=0.8\linewidth]{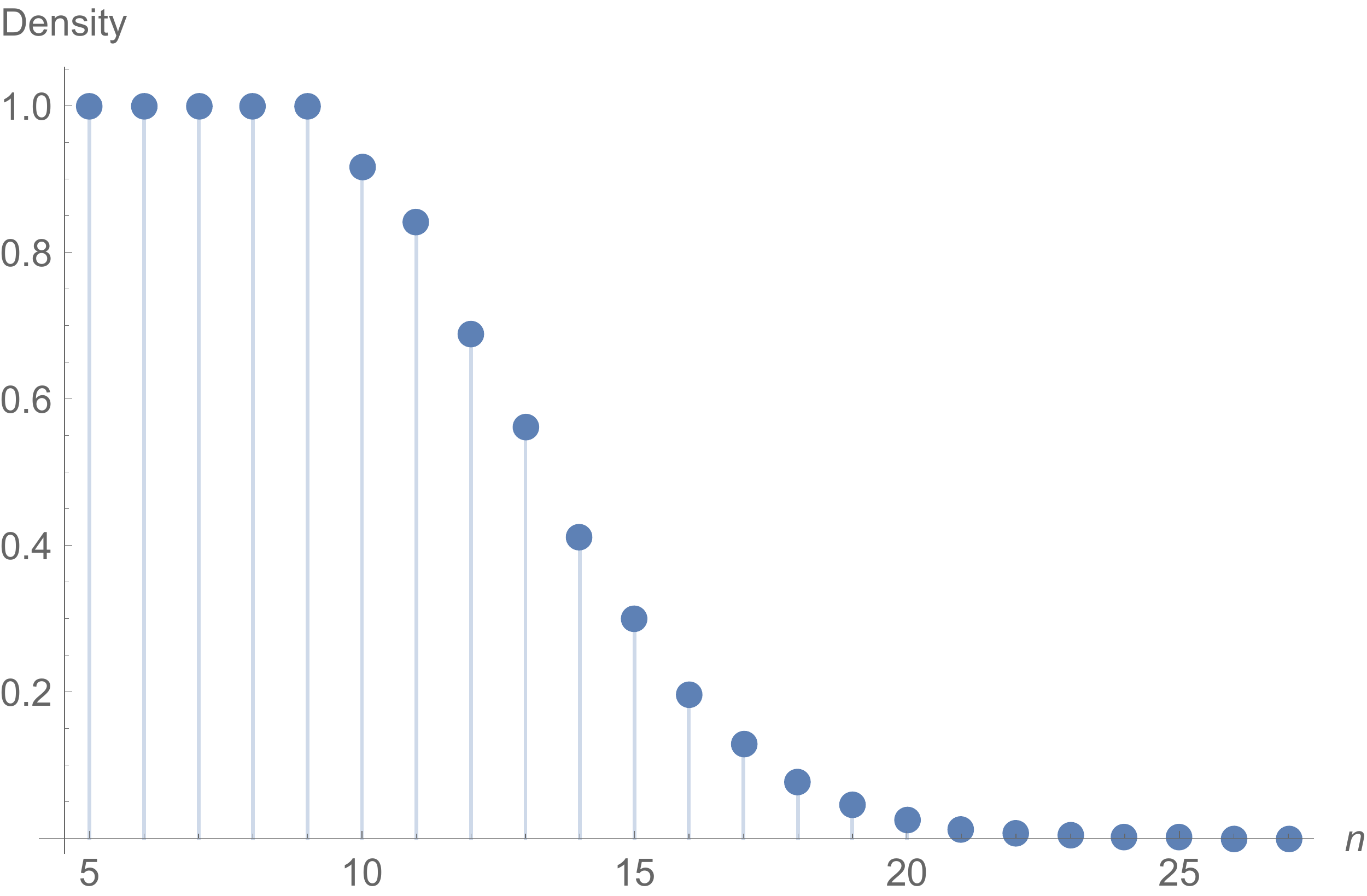}
	\caption{This plot shows the ratio of the number of non-zero entries to the total number of entries in a single KLT block. This shows that the generic behaviour is only reached for values of $n$ larger than $25$.}
	\label{densityKLT}
    \end{figure}

This problem is nothing but that of two copies of the problem solved for the complete intersection matrix. In other words, the number of non-zero intersections is given by 
\be
\# = (2\texttt{S}_{m-1})(2\texttt{S}_{\bar{m}-1})
\ee
with $\bar{m} = \lfloor (n-3)/2\rfloor +1$. The factor of $2$ in each part comes from the fact that reflections of each part must be counted separately, i.e. there is no cyclic symmetry for each part.

Using that the length of a row in a KLT block is given by $d$, the density is 
\be
\frac{4\texttt{S}_{m-1}\texttt{S}_{\bar{m}-1}}{(m-1)!(\bar{m}-1)!} = \left(\frac{\texttt{S}_{m-1}}{(m-1)!/2}\right) \left(\frac{\texttt{S}_{\bar{m}-1}}{(\bar{m}-1)!/2}\right).
\ee
The right hand side is written to show that the density of a KLT block is exactly the product of the densities of two intersection matrices, one with $n=m$ and another with $n=\bar{m}$.

The asymptotic behaviour can be read from \eqref{asym}  
\be
 \frac{\texttt{S}_{m-1}}{(m-1)!} \frac{\texttt{S}_{\bar{m}-1}}{(\bar{m}-1)!} \sim e^{-m(\log (m) - 2.76)-\bar{m}(\log (\bar{m}) - 2.76)}\sim  e^{-n(\log (n) - 3.46)}.
\ee
This shows that the KLT blocks are themselves sparse as $n$ is taken to be large. Once again, this conclusion is reached by computing the sparsity as discussed at the end of section 4. 

Repeating the same analysis as done for the full intersection matrix we consider the behaviour of the density of KLT blocks for finite values of $n$ and plot them in figure \ref{densityKLT}. Here we see that the point where we start seeing the generic behaviour is shifted to $n=25$.

\section{Discussions}

We hope that this note can motivate research on the combinatorial aspects of positive tropical Grassmannians embedded in ${\rm Trop}\, G(2,n)$ (for work on coverings of ${\rm Trop}\, G(2,n)$ by spaces of planar trees see \cite{Lauren1,Lauren2}).

Our main motivation has been the search for diagonal submatrices of the intersection matrix with the maximum possible rank. This maximum rank is what we called the diagonal degree of ${\rm Trop}\, G(2,n)$. 

A related problem is finding the maximum rank of diagonal submatrices when the two sets of orderings, ${\cal A}$ and ${\cal B}$, are required to be equal. While the maximum rank of such matrices is two for $n=5$ and coincides with the diagonal degree of ${\rm Trop}\, G(2,5)$, this is not the case already for $n=6$ where the maximum rank when the condition ${\cal A} = {\cal B}$ is imposed is only three\footnote{We thank N. Early for writing a Mathematica code which computed this number. In fact, the code was also able to find the maximum rank for $n=7$ which turns out to be six. This is to be compared to the lower bound we found for the diagonal degree of ${\rm Trop}\, G(2,7)$ which is fourteen.} while the diagonal degree of ${\rm Trop}\, G(2,6)$ is four. At this point we do not have reasons to believe that if the condition ${\cal A} = {\cal B}$ is imposed then the upper bound on the rank of diagonal submatrices, i.e. $(n-3)!$, will ever be reached for any $n>5$. Since our main interest is on physical applications, we have restricted the scope of our work to the diagonal degree of ${\rm Trop}\, G(2,n)$ as in Definition \ref{troprank}.     

In Example 2.1 we presented the results of an exhaustive search for $n=6$ where the intersection matrix is a $60\times 60$ matrix and found that the diagonal degree of ${\rm Trop}\, G(2,6)$ is $4$. While this search was easily done on a laptop, moving to $n=7$ requires new techniques as the number of possibilities increases very quickly. A preliminary search has led to a lower bound for the diagonal degree of ${\rm Trop}\, G(2,7)$ of $14$. Recall that in section 2 we proved that the upper bound is $(n-3)!$ which in this case turns out to be $24$.   

Of course, as shown in section 4, we do not expect to encounter the generic behavior of the diagonal degree of ${\rm Trop}\, G(2,n)$ until $n\sim 16$. This is clearly far beyond any exhaustive search can hope to reach and underscores the importance of finding purely combinatorial techniques to address this problem. The shift in behaviour around $n=16$ is reminiscent of a phase transition in physics. It is tempting to speculate that scattering amplitudes in theories constructed from the KLT procedure, such as Einstein gravity, exhibit a dramatic change in behaviour for large $n$.  

In \cite{CYY} it was suggested that using the closely related matrix of biadjoint amplitudes $m(\a,\b)$ used in section 2 to prove the upper bound on the diagonal degree, one could introduce the machinery of matroid and Ramsey theory into the problem. In other words, the vectors $\phi(\a)\in \mathbb{R}^{(n-3)!}$ have numerous properties which could given information about their linear dependence.  We leave these fascinating directions for future research.  

Instead, we end this work with some techniques that can be useful when computing intersection matrices and searching for maximal rank diagonal submatrices.

    \begin{figure}[h!]
	\hspace{-0.35in}
    \includegraphics[width=1.05\linewidth]{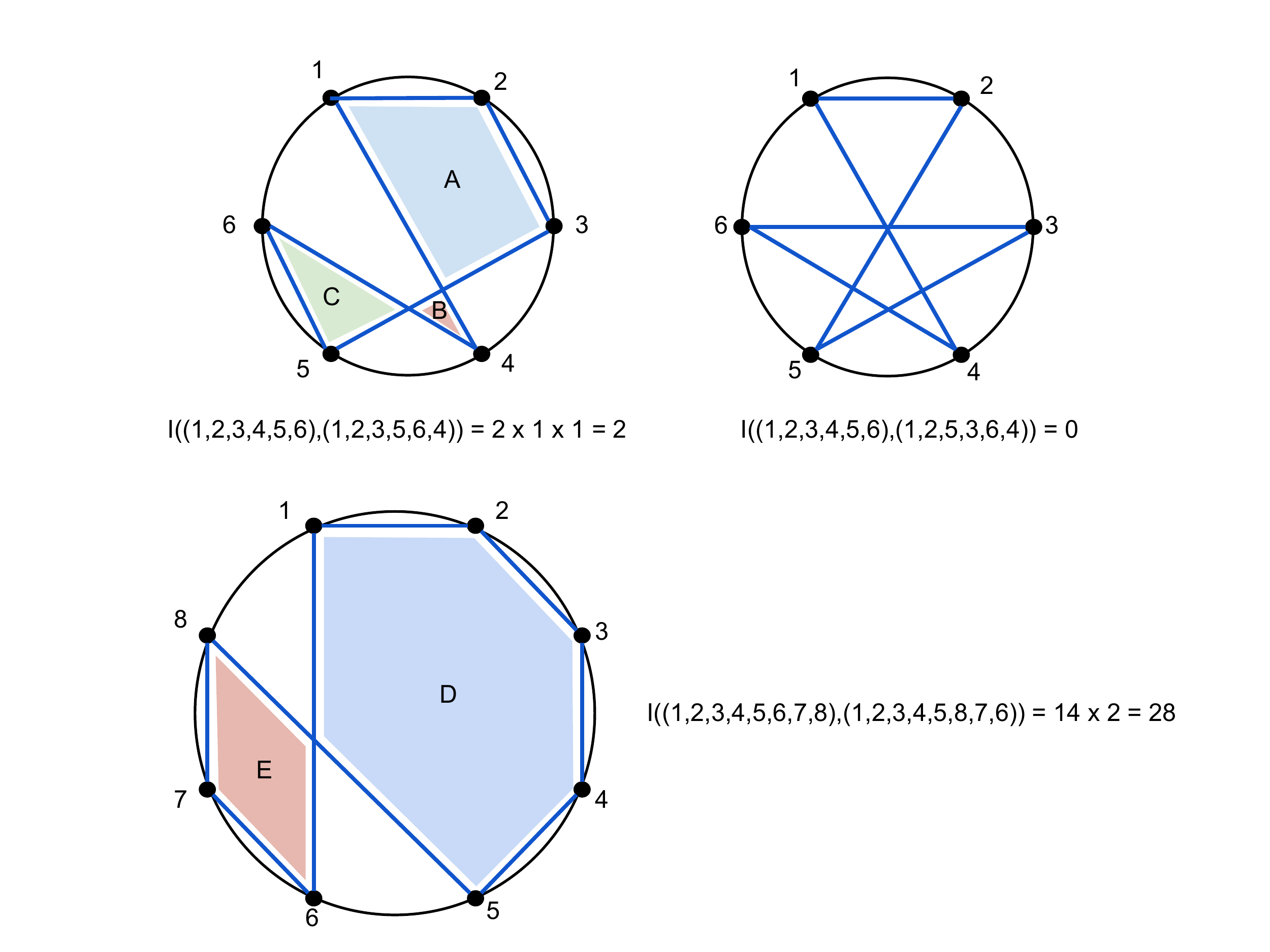}
	\caption{CHY diagrammatic computation of intersection numbers. Top left and bottom: Polygons are indicated with different colors. The contributions are: A and E are squares with $C_{2}=2$, B and C are triangles with $C_1=1$, and D is a hexagon with $C_4=14$. Top right: The polygons in the figure cannot be made to intersect only on vertices and therefore the intersection number vanishes.}
	\label{chyRules}
    \end{figure}

\subsection{CHY Diagrammatic Technique for Counting Intersections}

Consider elements of the intersection matrix in the row corresponding to the canonical order, i.e. of the form $I(\mathbb{I}_n , \a )$. In \cite{Cachazo:2013iea}, CHY introduced a diagrammatic technique which can be easily adapted for computing $I(\mathbb{I}_n , \a )$. The technique starts by drawing a circle and placing the labels in the order $\mathbb{I}_n = (1,2,\ldots ,n)$ on the boundary of the circle. Next take the order $\a = (a_1,a_2,\ldots, a_n)$ and draw straight lines from the location of $a_1$ on the boundary to the location of $a_2$ and so on until reaching $a_n$ and then ending back at $a_1$. If the resulting figure inside the circle can be drawn as a collection of polygons than only meet at vertices (allowing for continuous deformations of the straight lines), then the intersection number is given by 
\be
I(\mathbb{I}_n , \a ) = \prod_{i=1}^m C_{s_i-2}
\ee
where $m$ is the number of polygons, $s_i$ is the number of sides of the $i^{\rm th}$ polygon, and $C_r$ is the $r^{\rm th}$ Catalan number, otherwise the intersection number is zero.  

We illustrate this procedure with two examples taken from the computation of diagonal elements in Example 2.1 and a $n=8$ example. The diagrams are presented in figure \ref{chyRules}.

\subsection{Fast Algorithm for Detecting Zeroes}

In the search for diagonal submatrices of maximal rank of intersection matrices it is not necessary to compute the precise value of $I(\a,\b)$. The only information needed is whether an entry vanishes or not. Motivated by this let us define a less refined matrix: 

\begin{defn}\label{binary}
The binary intersection matrix $I^{\rm bin}(\a,\b)$ is a $(n-1)!/2\times (n-1)!/2$ matrix with entries
\be
I^{\rm bin}(\a,\b) = \left\{ \begin{array}{ll} 0 & \text{if } I(\a,\b) = 0, \\ 1  & \text{if } I(\a,\b) \neq 0. \end{array} \right. 
\ee
\end{defn}

There is a purely combinatorial way of computing $I^{\rm bin}(\a,\b)$ which is very useful in computer searches. 

The algorithm starts by computing the possible cherries a tree shared by both $\a$ and $\b$ can have. This is done by creating a list of all consecutive pairs of elements in both orderings and then finding the intersection. Each pair is to be treated as an un-ordered set. 

If there are no cherries in common then $I^{\rm bin}(\a,\b)=0$. Otherwise, select any of the cherries in common and prune it. This is done by taking the cherry, say given by the set $\{ i,j\}$, and removing one of the two labels, say $i$, from both $\a$ and $\b$ to produce two new orderings $\a_1 = \a \setminus \{i\}$ and $\b_1 = \b \setminus \{i\}$. 

The procedure is now repeated using $\a_1$ and $\b_1$ until either $I^{\rm bin}(\a_m,\b_m)=0$ for some $m$ or the length of $\a_m$ and $\b_m$ is four, in which case $I^{\rm bin}(\a_m,\b_m)=1$.

Let us illustrate the procedure with two examples. Consider first $\a = (1,2,3,4,5,6)$ and $\b =(1,3,6,5,2,4)$. The consecutive pairs obtained from each ordering are:
\be
\begin{array}{ccc}
    \a & \rightarrow & \{\{1,2\}, \{2,3\},\{3,4\},\{4,5\},\{5,6\},\{6,1\} \} , \\
    \b &  \rightarrow & \{ \{1,3\}, \{3,6\}, \{6,5\}, \{5,2\}, \{2,4\}, \{4,1\} \} .
\end{array}
\ee
Clearly the intersection is only a single pair, i.e. $\{6,5\}$. This provides our candidate cherry and it should be pruned. Pruning means removing the label $6$ to get $\a_1 = (1,2,3,4,5)$ and $\b_1=(1,3,5,2,4)$. Repeating the procedure shows that $\a_1$ and $\b_1$ do not share any pairs and therefore $I^{\rm bin}(\a,\b)=0$. 

The second example is given by $\a=(1,2,3,4,5,6)$ and $\b=(1,3,5,6,4,2)$. The common pairs are $\{ \{1,2\}, \{5,6\}\}$. Choose any of the cherries, say $\{1,2\}$ and prune it, i.e. remove label $1$ to get $\a_1=(2,3,4,5,6)$ and $\b_1=(3,5,6,4,2)$. Repeating the procedure gives the common pairs $\{ \{2,3\},\{5,6\}  \}$. Pruning $\{2,3\}$ by removing label $2$ gives $\a_2=(3,4,5,6)$ and $b_2=(3,5,6,4)$. Since the length of both $a_2$ and $b_2$ is equal to $4$ we stop and $I^{\rm bin}(\a,\b)=1$.

\subsection{Searching for a Maximal-Rank Permutation Submatrix in a Binary Matrix}

We end this work with a discussion regarding search algorithms that were implemented for the $n=6$ and $n=7$ cases. We are not aware of any publicly available algorithm so we include it here in the hopes that it can be useful or be improved in order to tackle higher $n$ cases.

Given a binary matrix $M$, the problem is to find the list of rows ${\cal R}$ and columns ${\cal C}$ with the largest possible number of elements such that the submatrix $M_{i,j}$ with $i\in {\cal R}$ and $j\in {\cal C}$ is a permutation matrix. Of course, by rearranging the lists ${\cal R}$ and ${\cal C}$ such a permutation matrix can always be turned into a diagonal matrix. 

The algorithm is the following: 

\begin{itemize}

\item Construct a list with the position of all $1$'s in $M$. Let's call it $\texttt{ONES}$.

\item Construct a graph $G$ with vertices labeled by elements of $\texttt{ONES}$ and include an edge between $(i,j)\in \texttt{ONES}$ and $(k,l)\in \texttt{ONES}$ if and only if neither $(k,j)$ nor $(i,l)$ belongs to $\texttt{ONES}$. Simply put, draw an edge if the rows $\{i,k\}$ and columns $\{j,l\}$ form a $2\times 2$ permutation matrix.

\item Find the largest clique in $G$ and call it $\texttt{CLIQUE}$.

\item Return ${\cal R} = \{ i : \exists \, j \, {\rm with}\, (i,j)\in \texttt{CLIQUE} \}$ and ${\cal C} = \{ j : \exists \, i \, {\rm with}\, (i,j)\in \texttt{CLIQUE} \}$.

\end{itemize}

It is easy to show that the algorithm returns the desired lists which generate the largest permutation submatrix of $M$. The advantage of this algorithm is that the search part has been turned into a well-known graph theoretic problem, ``find the largest clique", for which there are sophisticated implementations. Unfortunately, the clique problem is known to be NP-complete and running times for large graphs can greatly vary. Of course, here we have not used the intrinsic structure of the binary intersection matrix which can significantly reduce the search space. The simplest observation is the fact that any row is a permutation of any other and therefore we can select the first row as the first element of ${\cal R}$. This observation splits the problem into cases. For $n=7$ we have reduced the problem to that of searching for cliques in $17$ matrices of sizes of order $163\times 163$. A preliminary implementation which has randomly searched in roughly $20\%$ of the space has found diagonal submatrices of rank as large as $14$. While it is tempting to conjecture that $14$ is indeed the diagonal degree of ${\rm Trop}G(2,7)$, it could be that larger submatrices are very rare. An exhaustive search seems to be within reach and we leave this result for the future. Let us end this work by presenting one of the rank $14$ submatrices found so far:

\begin{align}
\nonumber {\cal A}\, & = &     
\{ (1, 2, 3, 4, 5, 6, 7), (1, 5, 6, 4, 3, 2, 7), (1, 3, 5, 4, 6, 2, 7), (1, 3, 4, 2, 6, 7, 5), (1, 2, 7, 6, 5, 4, 3), \\ \nonumber & & (1, 5, 3, 4, 2, 6, 7), (1, 2, 7, 3, 4, 5, 6), (1, 2, 7, 6, 4, 5, 3), (1, 3, 2, 4, 6, 5, 7), (1, 5, 3, 4, 7, 2, 6), \\ \nonumber & &  (1, 3, 2, 4, 6, 7, 5), (1, 2, 3, 4, 7, 5, 6), (1, 2, 6, 4, 3, 5, 7), (1, 5, 7, 3, 4, 2, 6) \}. \qquad \qquad \qquad \quad \\ \nonumber
{\cal B}\, & = & \{ (1, 4, 5, 2, 3, 6, 7), (1, 4, 6, 3, 7, 2, 5), (1, 3, 6, 2, 5, 4, 7), (1, 3, 6, 2, 5, 7, 4), (1, 3, 6, 5, 2, 7, 4),\\ \nonumber & & (1, 4, 2, 5, 3, 6, 7), (1, 4, 5, 2, 7, 3, 6), (1, 2, 5, 3, 6, 7, 4),
(1, 3, 6, 5, 2, 4, 7), (1, 4, 7, 2, 5, 3, 6), \\  & & (1, 4, 7, 6, 3, 2, 5), (1, 4, 7, 5, 2, 3, 6), (1, 2, 5, 7, 3, 6, 4), (1, 4, 2, 5, 7, 3, 6) \} .  \qquad \qquad \qquad \quad
\end{align}

\section*{Acknowledgements}

The author thanks A. Guevara, B. Sturmfels, S. Telen, and K. Yeats for discussions, D. Lang for help in implementing algorithms and especially N. Early for helpful discussions and many suggestions on the draft. We also thank B. Schr\"oter and L. Williams for comments on the draft. Finally, we thank S. Mizera and E. Yuan for numerous discussions on the structure of KLT blocks. This research was supported in part by a grant from the Gluskin Sheff/Onex Freeman Dyson Chair in Theoretical Physics and by Perimeter Institute. Research at Perimeter Institute is supported in part by the Government of Canada through the Department of Innovation, Science and Economic Development Canada and by the Province of Ontario through the Ministry of Colleges and Universities

\appendix

\section{Review of the CHY Formulation}

In this appendix we review the elements of the CHY formulation that are used in the proof of Proposition \ref{rank}. This review follows very closely section 2 of \cite{CYY}. 

There is a way of associating a rational function to every unrooted binary tree which is tightly connected to the tropical Grassmannian (or the BHV space of trees) and it is motivated by their physical application as Feynman diagrams in a cubic scalar quantum field theory. 

We start with a real $n\times n$ symmetric matrix, $s_{ab}$, satisfying the following properties
\be\label{con}
s_{aa} = 0\quad {\rm and}\quad \sum_{b=1}^n s_{ab}=0\quad \forall \; a\in \{1,2,\ldots ,n\}.
\ee
This space is $n(n-3)/2$ dimensional.

Let ${\cal T}$ be an unrooted binary tree with $n$ leaves and $E_{\cal T}$ be the set of edges connecting two trivalent vertices. Removing $e\in E_{\cal T}$ divides ${\cal T}$ into two disconnected graphs with a corresponding partition of the leaves into two sets $L_e\cup R_e=\{1,2,\ldots ,n\}$. Note that $|L_e|\geq 2$ and $|R_e|\geq 2$. The conditions in \eqref{con} imply that
\be
Q_e:= \sum_{a,b\in L_e}s_{ab} = \sum_{c,d\in R_e}s_{cd}
\ee
and therefore it is a quantity that can be associated with the edge $e$.

The rational function associated with ${\cal T}$ is then
\be
R({\cal T}) :=\prod_{e\in E_\Gamma}\frac{1}{Q_e}.
\ee

\begin{defn}\label{setCy}
Let $\Omega(\a)$ be the set of all unrooted binary trees with n leaves that a admit a planar embedding defined by $\a$, i.e. that belong to ${\rm Trop}^\a G(2,n)$.
\end{defn}

Now we are ready to give a formula for the object of interest.

\begin{defn}
A partial amplitude with orderings $\a$ and $\b$ is given by
\be\label{mdef}
m(\alpha , \beta) := (-1)^{w(\alpha,\beta)}\sum_{{\cal T}\in \Omega(\alpha)\bigcap\Omega(\beta)}R({\cal T}).
\ee
\end{defn}

In this formula the sum is over all trees that admit both a planar embedding defined by $\alpha$ and one defined by $\beta$. The overall sign is not relevant as we are only concerned by whether the object vanishes or not. We refer the interested reader to \cite{Cachazo:2013iea} for the definition of $w(\a,\b)$.

The CHY formulation of $m(\alpha , \beta)$ requires finding the critical points of 
\be
{\cal S}(x_1,x_2,\ldots ,x_n) :=\sum_{1\leq a<b\leq n}s_{ab}\, \log (x_a-x_b).
\ee
This is a Morse function on the moduli space of Riemann spheres with $n$ punctures. The number of critical points is given by the Euler characteristic of the space which is $(n-3)!$.
These $(n-3)!$ critical points are the solutions to what are known as the scattering equations \cite{Cachazo:2013gna,Cachazo:2013hca,Cachazo:2013iea}
\be\label{sceq}
\frac{\partial {\cal S}}{\partial x_a} = \sum_{b=1,b\neq a}\frac{s_{ab}}{x_a-x_b}=0\quad \forall\; a\in \{1,2,\ldots,n\}.
\ee
Let's denote the $(n-3)!$ solutions as $x_a^{I}$. In general the solutions are complex but when the $s_{ab}$'s are chosen in what is known as the positive region all solutions are real \cite{Cachazo:2016ror}. Given any ordering $\alpha$ one constructs a vector $\phi(\alpha)\in\mathbb{R}^{(n-3)!}$ whose components are given by
\be\label{defPhi}
\phi(\alpha)_I := \frac{K_I}{(x_{\alpha_1}^I-x_{\alpha_2}^I)(x_{\alpha_2}^I-x_{\alpha_3}^I)\cdots (x_{\alpha_n}^I-x_{\alpha_1}^I)},
\ee
where $K_I$ is a function obtained from second derivatives of ${\cal S}$ and it is invariant under permutations of labels and hence $\alpha$ independent. Therefore $K_I$ is not relevant to our discussion and we refer the reader to \cite{Cachazo:2013iea} for details.  

Finally, partial amplitudes are computed as
\be\label{weq}
m(\alpha , \beta) = \sum_{I=1}^{(n-3)!}\phi(\alpha)_I\,\phi(\beta)_I.
\ee
We also use the notation $\phi(\alpha)\cdot\phi(\beta)$ for the inner product in \eqref{weq} in the main text. This is the formula used in the proof of Proposition \ref{rank} in section 2. This formula was originally proposed by He, Yuan, and the author \cite{Cachazo:2013gna,Cachazo:2013hca,Cachazo:2013iea} and later proven by Dolan and Goddard in \cite{Dolan:2013isa}.

\section{Asymptotic Behavior of Super Catalan Numbers}

In section 4 we used the asymptotic behaviour of the super Catalan numbers. The reference cited in the main text shows the result but not the proof and since we could not find a simple derivation in the literature we present one here for completeness. Let us start with the recursion relation $\texttt{S}_{r}$ is known to satisfy,
\be\label{oneR}
\texttt{S}_{r} = \frac{3(2r-3)\texttt{S}_{r-1}-(r-3)\texttt{S}_{r-2}}{r}
\ee
with $\texttt{S}_{1}=\texttt{S}_{2}=1$ \cite{wolfram}.

Now let us recall the recursion relation the Legendre polynomials satisfy:
\be\label{twoR}
P_{l+1}(x) =\frac{(2l+1)x\,P_l(x)-l\,P_{l-1}(x)}{l+1}.
\ee
The similarities between \eqref{oneR} and \eqref{twoR} motivate a connection between $\texttt{S}_{n}$ and $P_l(x)$ evaluated at $x=3$. 

In fact, it is not difficult to show that a linear combination of Legendre polynomials satisfies the recursion \eqref{oneR} and the boundary conditions \cite{wolfram}
\be\label{aqe}
\texttt{S}_{r}=\frac{3P_{r-1}(3)-P_{r-2}(3)}{4r}.
\ee
The asymptotic behaviour of Legendre polynomials $P_l(x)$ as $l$ becomes large and $x>1$ is known to be
\be\label{asymP}
P_l(x) = \frac{1}{\sqrt{2\pi l y(x)}}\frac{(1+y(x))^{(l+1)/2}}{(1-y(x))^{l/2}}+{\cal O}(l^{-1}),
\ee
where
$$ y(x) = \sqrt{1-1/x^2}. $$
Using $x=3$ gives $y=\sqrt{8}/3$. Combining this with \eqref{asymP} and \eqref{aqe} gives the asymptotic behavior used in section 4,
\be
\log(\texttt{S}_{n-1})\sim n\, \log\left(3+\sqrt{8}\right)-\frac{3}{2}\log(n) + {\cal O}(n^0). 
\ee

\bibliographystyle{JHEP}
\bibliography{references}

\end{document}